\documentclass{amsart}
\usepackage[latin1]{inputenc}
\usepackage{amssymb}
\usepackage{latexsym}
\usepackage{longtable}
\usepackage{graphicx}
\DeclareMathAlphabet{\mathfr}{U}{euf}{m}{n}
\newtheorem{theorem}{Theorem}[section]

\newtheorem{rem}{Remark}[section]

\newtheorem{lemma}{Lemma}[section]

\newtheorem{proposition}{Proposition}[section]

\newcommand{\Q}{\mathbb Q}

\newcommand{\Gal}{\mathrm{Gal}}

\newcommand{\Z}{\mathbb Z}

\newcommand{\C}{\mathbb C}

\newcommand{\cO}{\mathcal{O}}
\newcommand{\ra}{{\rightarrow}}
\newcommand{\cL}{\mathcal{L}}
\newcommand{\End}{\operatorname{End}}

\newcommand{\Aut}{\operatorname{Aut}}
\newcommand{\re}{\operatorname{Re}}
\newcommand{\II}{\operatorname{Im}}

\newcommand{\Res}{\operatorname{Res}}

\newcommand{\Jac}{\operatorname{Jac}}
\newcommand{\NS}{\operatorname{NS}}
\newcommand{\im}{\mathrm{Im}}

\newcommand{\qbar}{\overline{\Q}}

\newcommand{\ds}{\displaystyle}

\begin{document}

\title{Genus two curves with quaternionic multiplication and modular jacobian}
\author{Josep Gonz\'alez and Jordi Gu\`ardia}
\date{\today}

\address{Escola  Polit\`ecnica  Superior d'Engenyeria  de Vilanova i la Geltr\' u,
Avda Victor Balaguer s/n,  08800 Vilanova i la Geltr\'u, Spain}
\email{josepg@ma4.upc.edu, guardia@ma4.upc.edu}

\thanks{The  authors are partially supported by  MTM2006-15038-C02-02.}

\maketitle
\begin{abstract}
We describe a method to determine all the isomorphism classes of
principal polarizations of the modular abelian surfaces $A_f$  with
quaternionic multiplication attached to a normalized newform $f$
without complex multiplication. We include an example of $A_f$ with quaternionic multiplication for which we find numerically a curve $C$
whose Jacobian is $A_f$ up to numerical approximation, and we prove that it has quaternionic multiplication and is  isogenous to  $A_f$.

\end{abstract}

\section{Introduction}

During the last years, abelian surfaces have emerged in arithmetic geometry in different contexts,  illustrating interesting phenomena and providing examples  for practical applications. While the general theory of abelian surfaces is well-known, explicit methods are quite recent and not completely developed.
For an irreducible abelian surface over a number field, the problem of determining  explicitly all the genus two curves whose Jacobian is isomorphic to the given surface is not solved in general.

We began the study of this problem for irreducible principally polarized abelian surfaces in \cite{GGG}.
Afterwards, in \cite{GGR} we developed the theoretical results related to the problem for irreducible polarized abelian surfaces, considering specifically non-principal polarizations. In both cases, we applied our ideas to modular abelian surfaces, since the apparition of new algorithms for the description of the Jacobians of modular curves makes possible the generation of explicit examples. Unfortunately, the numerical nature of this algorithms drives to  only numerically tested examples.

We now consider the   interesting case of abelian surfaces with quaternionic multiplication, whose rich endomorphism algebra allows the combined application of further well-known techniques. After describing these techniques, we provide an example of a modular abelian surface presented as the Jacobian of a hyperelliptic curve for which,  although numerically determined, we can determine its endomorphism ring and prove that it is isogenous to the given surface.

\section{General facts}

Let   $f=\sum_{n>0}a_n q^n$  be a normalized newform of $S_2(\Gamma_0(N))$ without complex multiplication  such  that the abelian variety $A_f/\Q$ attached by Shimura in \cite{shimura71} has
dimension~$2$.
The $\Q$-algebra $\End^0(A_f):=\End_{\overline
\Q}(A_f)\otimes \Q$  can only  be the real quadratic field generated
by the Fourier coefficients $a_n$, the matrix algebra $\operatorname{M}_2 (\Q)$ or an
indefinite quaternion algebra $(\frac{a,b}{\Q})$. We restrict
ourselves to the most interesting case that $A_f$ has quaternionic multiplication. Indeed, we will see that if there is a
principal polarization not defined over $\Q$, $A_f$ can be presented
as the jacobian of a genus two-curve $C$ defined over an imaginary
quadratic field which does not admit descent over $\Q$. Then, there
are two possibilities: either $C$ is isomorphic to its complex
conjugate $\overline C$ and, thus, with rational Igusa invariants or
the isomorphism classes of the canonical principal  polarizations of
both curves are different.

 We denote by $E=\Q(\sqrt m)$ the real
quadratic number field generated by $\{a_n, n>0\}$, where $m$ is a
square-free integer. It is well known that every Hecke operator
provides an endomorphism of $A_f$ for which $f$ is an eigenvector
and that we identify with the eigenvalue of $f$. So we have a natural injection $\Z[ \{
a_n\}]\hookrightarrow \End_\Q(A_f)\subset \End_{\Q}^0(A_f)\simeq E$ which let us interpret the coefficients $a_n$ as endomorphisms of $A_f$.

An inner twist of $f$  is a Dirichlet character $\chi$ of conductor
dividing $N$ such that ${}^\sigma f=f\otimes \chi$ for some
nontrivial $\sigma\in\Gal (\overline \Q/\Q)$, i.e. $a_p=\chi
(p)\,{}^{\sigma}a_p$ for all primes $p\nmid N$ . Since $f$ has
trivial Nebentypus,  $\chi$ must be quadratic (cf. \cite{ribet80}); we will write $K=\overline{\Q}^{\ker \chi}$. By Proposition 8 of
\cite{shimura73} there is $u_\chi\in \End_{K}^0 (A_f)$ sending $f$
and ${}^{\sigma}f$ to $g(\chi) {}^{\sigma}f$ and $g(\chi) f$
respectively, where $g(\chi)$ denotes the Gauss sum relative to the
conductor of $\chi$. Of course, $u_\chi^2=\operatorname{disc } K$.
From \cite{ribet80}, we have the following characterization.

\begin{proposition}\label{algebra} The abelian surface $A_f$ has quaternionic
multiplication if and only if there is an inner twist $\chi$ such
that its conductor is not  a norm of the number field $E$. In this case,
the quadratic number field $K=\Q(\sqrt \delta\,)=\overline{\Q}^{\ker
\chi}$ is imaginary and $\End^0(A_f)=\End_K^0(A_f)=\Q(u_\chi,
T_p)\simeq\left ( \frac{\delta, m}{\Q}\right)\,,$ where $T_p$ is the
Hecke operator at a prime $p$  such that $a_p\notin\Z$.

\end{proposition}




Let $k$ be a subfield of a  fixed algebraic closure $\overline \Q $
of $\Q $, let  $G_k$ be the absolute Galois group $\Gal (\overline
\Q /k)$ and let $A$  be an abelian variety defined over $k$. The
class of an invertible sheaf $\cL$ (not necessarily defined over $k$)
in the N\'eron-Severi group $\NS (A_{\overline \Q})$ is defined over
$k$ if it is invariant by $G_k$. Every invertible sheaf $\cL \in \NS
(A_{\overline \Q})^{G_k}$ defines a morphism $\varphi _{\cL }: A \ra
\hat A$ over $k$ given by $\varphi _{\cL }(P)=\tau _P^*(\cL )\otimes
\cL ^{-1}$, where $\tau_P$ denotes the translation by $P$. This
morphism is an isogeny if and only if $\cL $ is nondegenerate. A polarization on
$A$ defined over $k$ is the class of algebraic equivalence of an
ample invertible sheaf $\cL \in \NS (A_{\overline \Q})^{G_k}$.
Equivalently, a polarization on $A$ over $k$ is an isogeny $\lambda
:A \ra \hat A$ defined over $k$ such that $\lambda \otimes \overline
k = \varphi _{\cL }$ for some ample line bundle $\cL $ on
$A_{\overline \Q}$. The polarization $\cL$ is called  principal when
the degree of $\varphi_\cL$ is $1$.

We proceed to present the theoretical that we will use for our computations.

\begin{proposition} [cf.  \cite{weil57}]
Let $A/k$ be an abelian surface with a principal polarization $\cL $
defined over  $k$. If $A$ is simple over $\overline \Q$, then $(A,
\cL ) \stackrel{k}\simeq (\Jac (C), \cL (\Theta _C))$, where $C/k$
is a smooth curve of genus two and $\cL (\Theta _C)$ denotes the
canonical principal polarization for $\Jac (C)$.

\end{proposition}

A nondegenerate invertible sheaf $\cL $ on $A$ defined over $k$
induces an anti-involution on the algebra of endomorphisms
$$
*:\Q \otimes \End _{k}A\stackrel{\sim }{\ra }\Q \otimes \End
_{k}A, \quad t\mapsto {\varphi _{\cL }}^{-1}\cdot \hat t\cdot
\varphi _{\cL }.
$$
Let $\End _k^sA=\{ \beta \in \End _k A, \beta ^*=\beta \} $ denote
the subgroup of symmetric endomorphisms and $\End _{k +}^sA$ be the
set of positive symmetric endomorphisms of~$A$.

\begin{theorem}[Proposition 2.1 and Theorem 2.3 of \cite{GGR}]
\label{pol}

Let $A/k$ be an abelian variety and let $\cL \in \NS (A_{\qbar
})^{G_k}$ be nondegenerate. Then,

\begin{itemize}

\item[(i)] For any endomorphism $t\in \End_k^s A$, there exists a unique
$\cL ^{(t)}\in \NS (A_{\qbar })^{G_k}$ such that $\varphi _{\cL
^{(t)}}= \varphi _{\cL }\cdot t$. More precisely, if  $E$ and $E_t$
denote the alternating Riemann forms attached to $\cL$ and
$\cL^{(t)}$ respectively. Then
$$
E_t(x,y)=E(x, t y)=E(t x, y)\,.
$$
Moreover, if $t$ is a totally positive element, then $\cL$ is a
polarization if and only if $\cL^{(t)}$ is.

\item[(ii)] Assume $A$ is principally polarized over $k$.  For any choice of a principal polarization $\cL _0$ on $A$ defined
over $k$, there is an isomorphism of groups

$$
\begin{matrix}
\epsilon : &\NS (A_{\qbar })^{G_k}&\stackrel{\sim }{\ra }&\End _k^sA\\
         &   \cL   &\mapsto       &  \varphi _{\cL _0}^{-1}\cdot
          \varphi _{\cL }
\end{matrix}
$$
such that $\cL \in \NS (A_{\qbar })^{G_k}$ is a polarization if and
only if $\epsilon (\cL )\in \End _{k +}^sA$ and it is principal if
and only if $\epsilon (\cL )\in \Aut _{k +}^sA$. Moreover,
$\epsilon^{-1}(t)=\cL_0^{(t)}$. Moreover, the set all of
$k$-isomorphism classes of principal polarizations on $A$ defined
over $k$ is the set  $\epsilon^{-1}(\Aut _{k +}^sA/\sim)$, where
$\beta _1\sim \beta _2$ means that  $\beta _1= \beta ^*\beta _2
\beta $ for some $\beta \in \Aut _{k}A$, is the set of  .

\item[(iii)] Let $\cL$ be a polarization on $A$ over $k$ of degree
$d\geq 1$. Then, $A$ is principally polarizable over a number field
$k$ if and only if there exists $\gamma\in   \End_{k+}^s(A)$
satisfying $\deg \,\gamma =d^2$ and $ \cL ^{(\gamma^{-1})}\in \NS
(A_k)^{G_k}$.

\end{itemize}

\end{theorem}
\noindent  {\bf Proof.}  Part (i)    and part (ii) can be found in
Proposition 2.1 and Theorem 2.3 of \cite{GGR} respectively. Part
(iii) follows from part (ii) by using  the same arguments as in
Corollary 2.12 of \cite{GGR}  and taking into account that  $\deg
\cL=d$ is equivalent to  $\deg \varphi_{\cL}=d^2$. \hfill $\Box$

In the particular case that $A$ is a surface and $\End^0(A)$ is an
indefinite quaternion algebra, we have the following results.

\begin{theorem}[Theorem 1.1 of \cite{Rotger}]
\label{rotger}
Assume that $\End(A)$ is a maximal  order  of  an indefinite quaternion
algebra of discriminant $D$. Then, $A$ is principally polarized and
the number $\pi(A)$ of isomorphism classes of principal polarizations of $A$,
 is
$$
\pi(A)=\left\{
\begin{array}{ll}
\displaystyle\frac{h(-4\,D)+ h(-D)}{2} & \text{ if $D\equiv-1\pmod 4$}\\[ 3 pt]
\displaystyle\frac{h(-4\,D)}{2} & \text{ otherwhise.}\end{array}\right.
$$

\end{theorem}

\begin{proposition}[Lemma 4.3 of \cite{BFGR}]
\label{pp}
Let $A/\Q $ be an abelian surface such that $\End_\Q^0(A)$ is the
real quadratic  number field $\Q(\sqrt m )$ and $\End (A)$ is a
maximal order in a quaternion algebra of discriminant $D$. Then, $A$
admits a polarization of degree $d
>0$ defined over $\Q$ if and only if $\End_k^0(A)\simeq (\frac {-D d, m}{\Q
})$. In particular,~$A$~is the jacobian of a genus 2 curve defined
over $\Q$  if and only if the algebra $(\frac {-D , m}{\Q })$
ramifies exactly at the primes dividing $D$.
\end{proposition}

\section{Determination of curves with  jacobian isomorphic to $A_f$}

Let us write $E=\Q (\sqrt m )$, $K=\Q(\sqrt {\delta})$ with $m>0,
\delta<0$  square-free integers. We denote by  $i\in \End_\Q^0(A_f)$
the (fractional) Hecke operator such that $i^2=m$.  We take
$j\in\End^0_K(A_f)$ as $u_\chi$  if the discriminant of $K$ is
$\equiv 1\pmod 4$ and as or $u_\chi/4$ otherwise. Then $j^2=\delta$
and we set $k:=i\cdot j=-j\cdot i$. Let us denote  $\cO$  the order
$\End (A_f)$.

In order to apply theorem \ref{pol}, we have to determine first the order  $\cO$ in  $\Q(i,j)$. Once we know $\End
^s(A_f)$ and $\End_{+}^s(A_f)$, we  are able to compute the number
of isomorphism classes of principal polarizations and
 we can determine
hyperelliptic equations for each one of these polarizations by
applying the procedure described in \cite{GGG}.

\subsection{The ring of endomorphims of $A_f$}
Let $A/k$ be an abelian variety of dimension $n$ over a subfield $k$
of $\C$, $\{\omega_1,\cdots,\omega_n\}$ and $\{c_1,\cdots, c_{2n}\}$
be arbitraries  bases of
 $H^0(A,\Omega_{A/k}^1)$ and
 of $H_1(A,\Z)$ respectively. We can take as period lattice for $A$:
$$\Lambda=\left\{ \left(\left.\int_{c}\omega_1,\cdots
,\int_{c}\omega_n\right)\right| c \in H_1(A,\Z)\right\}=
\left\langle\left(\int_{c_i}\omega_1,\cdots,\int_{c_i}\omega_n\right), i\leq 2 n
\right\rangle.$$ The period  matrix $\Omega$ of $A$ with respect these
bases $\omega_1,\cdots,\omega_n$ and $c_1,\cdots,c_{2n}$ is
$$
\Omega=\left(\int_{c_i} \omega_j\right)_{1\leq i \leq 2n\,,\,1\leq j
\leq n}=\left(
\begin{array}{ccc}
\int_{c_1} \omega_1&\dots    &\int_{c_{2n}} \omega_1\\
  \vdots  & \vdots  & \vdots \\
  \int_{c_{1}} \omega_n&\dots    &\int_{c_{2n}} \omega_n
\end{array}
\right )
 $$
 An endomorphism $\phi:A\rightarrow A$ defined over $k$ induces an endomorphism $\phi^*$ in the $k$-vector
 space $H^0(A,\Omega_{A/k}^1)$ and an endomorphism $\phi_*$ in the
 $\Z$-module $H_1(A,\Z)$. The action of both endomorphisms is
 related by the equality
 \begin{equation}\label{action}
 \int_{\phi_*(c)}\omega=\int_c \phi^*(\omega)\quad
 \text{for all } \omega\in H^0(A,\Omega^1) \quad \text{and}\quad
 c\in H_1(A,\Z)\,.
 \end{equation}

 Let us denote by $T\in \operatorname{M}_n(k)$ and $M\in
 \operatorname{M}_{2n}(\Z)$ be the matrices of $\phi^*$ and $\phi_*$ with respect to
 the chosen bases. It follows from (\ref{action}) that
 \begin{equation}\label{endo1}
{}^tT\cdot \Omega=\Omega\cdot M \,\,.
 \end{equation}
Therefore
\begin{equation}\label{endo2}
M= \left(
\begin{array}{c}
\re \Omega\\
\II \Omega
\end{array} \right )^{-1}.
\left(
\begin{array}{c}
\re ({}^tT\cdot \Omega)\\
\II ({}^tT\cdot \Omega)
\end{array} \right )\,.
\end{equation}

 A $\C$-linear endomorphism $\psi$ on $ H^0(A,\Omega^1)$ comes from
some $\phi\in \Q\otimes \End(A)$, i.e., $\psi=\phi^*$, if and only
if $\psi (\Lambda)$ is contained in  $\Q\otimes \Lambda$, where
$\psi (\Lambda)$ is
$$\{ (\int_{c}\psi(\omega_1),\cdots,\int_{c}\psi(\omega_n))\mid c \in H_1(A,\Z)\}=
\langle
(\int_{c_i}\psi(\omega_1),\cdots,\int_{c_i}\psi(\omega_n))\,,\,\,
1\leq i \leq 2n \rangle\,.
$$

As an immediate consequence of the equalities (\ref{endo1}) and
(\ref{endo2}), we obtain the next criterion.

\begin{lemma}\label{criterion}

If $T$ is the matrix of a $\C$-linear endomorphism $\psi$ of $H^0(A,\Omega^1)$ with respect to the basis $\{\omega_1,
\cdots,\omega_n \}$, then $\psi=\phi^{*}$ for some $\phi\in \Q
\otimes \End (A)$  if and only the matrix
$$M_T:=
\left(
\begin{array}{c}
\re \Omega\\
\II \Omega
\end{array} \right )^{-1}.
\left(
\begin{array}{c}
\re ({}^tT\cdot \Omega)\\
\II ({}^tT\cdot \Omega)
\end{array} \right )
$$
lies in $ \operatorname{M}_{2n}\,(\Q)$ and, moreover,
$${}^tT.\Omega=\Omega.M_T \,.$$
When this is the case, $\phi\in  \End (A)$ if and only if $M_T\in
\operatorname{M}_{2n}\,(\Z)$.
\end{lemma}

In our case $A=A_f$, Proposition \ref{algebra} gives an explicit
description of the action of $\End^0(A_f)=\Q(i, j)$ on
$H^0(A_f,\Omega^1)$. Then,  to determine the order $\cO$ in
$\Q(i,j)$ we will only use fractional endomorphisms and, thus, we
will only need to check the condition $M_T\in
\operatorname{M}_{2n}\,(\Z)$ for a finite number of matrices $T\in
\operatorname{M}_{n}\,(K)$.

\subsection{Principal polarizations} We chose  a
basis $f_1,f_2$  of the vector space $\langle f
{}^{\sigma},f\rangle$ with rational Fourier coefficients. We take
$\omega_i=f_i(q) \,dq/q$, $i\leq 2$, as a basis of
$H^0(A_f,\Omega^1_{A_f/\Q})$. We compute with {\em Magma}  (\cite{magma}) the data necessary to describe the subvariety $A_f$: the modular symbols giving a basis of pathes for the homology, the period matrix of $\omega_1, \omega_2$ and the matrix $M_E$ of the alternating Riemann form $E$ (with
respect to this basis) attached to the canonical polarization $\cL$
of the jacobian of $X_0(N)$.

Since $\cL$ is defined over $\Q$, the Rosati involution  $*$
attached to $\cL$ acts on $\End_\Q^0 (A_f)$ as the complex
conjugation and hence
$$
 i^*=i\,, \quad  (j^*)^2=\delta\,,\quad i\cdot j^*=-j^*\cdot i\,,\quad  k^*=j^*\cdot i^*\,.
 $$
Therefore  $j^*=\pm j$ and $k^*= \mp k$ and $\End^s (A_f)$ is a free
$\Z$-module of rank $3$. The action of $*$ on $\End^0(A_f)$ viewed
as a subfield of $ \operatorname{M}_4(\Q) \simeq \End (\Q\otimes
H_1(A_f,\Z))$ is the following
$$M\mapsto  ( M_E\cdot M\cdot
M_E^{-1})^t\,.$$ So, once we have determined a $\Z$-basis of the
order $\cO$,  we have the corresponding matrices $M_i\in
\operatorname{M}_4(\Z)$ acting on $H^1(A_f,\Z)$ with respect to the
basis chosen for this $\Z$-module. Then $\End^s(A_f)$ is given by
the matrices
$$ \{ M\in \Z M_1\oplus \cdots \oplus\Z M_4:  M^t\cdot M_E=M_E \cdot
M\}\,,$$ and $\End_+^s(A_f)$ is the subset corresponding to the
matrices $M$ with  positive eigenvalues.

When the polarization $\cL$ is of type $(1,d)$ with $d>1$, using part (iii)
of Theorem \ref{pol} we can check whether $A_f$  admits
principal polarizations and
compute all isomorphism classes of them.

\subsection{Rational equations}

Once we know that an irreducible abelian surface $A_f$ is principally polarizable, we can realize it as the Jacobian variety of a hyperelliptic curve over $K$. We find a $K$-rational model $Y^2=F(X)$  of this curve following the procedure described in \cite{GGG}, with two slight corrections. First, the computation of the leading coefficient of the polynomial $F$ can be speeded up using directly Thomae's formula (cf. theorem 11.1 in \cite{guardia2002}).  Secondly, the polynomial $F(X)$ is determined up to a fourth root of unity in $K$ (not only a sign, as stated in \cite{GGG}): the explicit determination of this unit is performed anyway as in {\em{loc. cit.}}

\section{Example:  $A_f$ with only one principal polarization}

Let us consider the abelian surface $A_f$ attached to the normalized newform in $S_2(\Gamma_0(243))$
$$
f=q + \sqrt{6}q^2 + 4q^4 - \sqrt{6}q^5 + 2q^7 + 2\sqrt{6}q^8 - 6q^{10} + \sqrt{6}q^{11} - q^{13} +
    2\sqrt{6}q^{14} + 4q^{16}+\cdots
$$
We take the  basis of regular differentials over $\Q$ given by
$$
\omega_1=\frac12(f+{}^\sigma f)\frac{dq}{q},\qquad \omega_2=\frac1{2\sqrt{6}}(f-{}^\sigma f)\frac{dq}{q}.
$$
The polarization $\cL$ on $A_f$ induced by the canonical
polarization on the jacobian of $X_0(243)$ has type $[1,2]$. The
period matrix of $A_f$ with respect to a symplectic basis $\{c_1,
c_2, c_3, c_4\}$ of $H_1(A_f,\Z)$ and $\omega_1,\omega_2$ is
$\Omega=(\Omega_1|\Omega_2)$ with
$$
\begin{array}{c}
\Omega_1=
\left(
\begin{array}{cc}
12.3724\dots + 21.4297\dots I & 13.5178\dots  + 7.8045\dots I \\
 -2.253\dots  + 3.9023\dots I &-12.3724\dots  + 7.1432\dots I
\end{array}
\right),
\\ \\
\Omega_2=
\left(
\begin{array}{cc}
     -12.3724\dots  + 21.4297\dots I  &  - 15.609\dots I\\
 2.253\dots  + 3.9023\dots I&  - 14.2865\dots I
\end{array}
\right),
\end{array}
$$
where $I=\sqrt{-1}$. Over $K=\Q(\sqrt{-3})$ the endomorphism algebra
of $A_f$ is $\End_K^0(A_f)= \left ( \frac{6, -3}{\Q}\right)= \Q(
i,j)$, with $i^2=6, j^2=-3$. Applying formula (\ref{endo2}) we
determine the matrices giving the action of  the quaternions $i,j $
on $H_1(A_f,\Z)$. We  find:
$$
i:
\left(\begin{array}{rrrr}
 0 &-2 & 0 &-2\\
-1 & 0 & 1 & 0\\
 0 & 4 & 0 &-2\\
-2 & 0 &-1 & 0
\end{array}\right), \qquad
j:\left(\begin{array}{rrrr}
-1&  0& -2 &\phantom{-} 0\\
 0& -1&  0 & 2\\
 2&  0&  1 & 0\\
 0& -2&  0 & 1\\
\end{array}\right).
$$
The criterion given in Lemma \ref{criterion} shows that
$\frac{1}{2}i+\frac{1}{6} k,\frac{1}{2}+\frac{1}{2} j, \frac{1}{3}
k$ are also endomorphisms of $A_f$. Since $\cO=\langle1,
\frac{1}{2}i+\frac{1}{6} k,\frac{1}{2}+\frac{1}{2} j, \frac{1}{3}
k\rangle$ is a maximal order in $\left ( \frac{6, -3}{\Q}\right)$,
it must be  $\End_K(A_f)=\cO$. Hence, by Theorem \ref{rotger} there
exists only one isomorphism class of principal polarizations on
$A_f$; it is defined over $K$ but not over $\Q$ by Proposition \ref{pp}. A
principal polarization must be of the form
 $\cL^{(\gamma^{-1})}$ for
a totally positive symmetric endomorphism $\gamma\in\cO$ with
reduced norm 2. Indeed, we take
$$
\gamma=2+\frac12i-\frac16k:\left(\begin{array}{rrrr}
 2& \phantom{-} 0&  \phantom{-}0& -2\\
 0&  2&  1&  0\\
 0&  2&  2&  0\\
-1&  0&  0&  2\\
\end{array}\right).
$$
The Riemann form attached to $\cL^{(\gamma^{-1})}$  is
$E_{\gamma^{-1}}(x,y)=E(\gamma^{-1}x,y)$ where $E$ is the Riemann
form attached to $\cL$. Thus, in the basis $\{c_1, c_2, c_3, c_4\}$
of $H_1(A_f,\Z)$ it is given by the matrix
$$E_\gamma: \left(\begin{array}{rrrr}0 & 0 & \phantom{-}1 &\phantom{-}0 \\ 0 & 0 & 0 & 2 \\ -1 & 0 & 0 & 0\\ 0 & -2 & 0 & 0
\end{array}\right)\gamma^{-1}=
\left(\begin{array}{rrrr}
 0& -1&  \phantom{-}1 &\phantom{-} 1\\
 1&  0&  1 & 2\\
-1& -1&  0 & 0\\
-1& -2&  0 & 0\\
\end{array}\right).
$$
With respect to $E_{\gamma^{-1}}$,   a symplectic basis of is
$(c_1', c_2', c_3', c_4')= {}^t M. (c_1, c_2, c_3, c_4)$, with
$$
M= \left(\begin{array}{rrrr}  0  &1 & \phantom{-}1&  0\\ 0 &-1&  0&
0\\-1  &0 & 0&  1\\ 0& -1&  0& -1
\end{array}\right),
$$
so that a period matrix for $(A_f,\cL^{(\gamma^{-1})})$ is
$\Omega_\gamma:=\Omega. M$. We can finally apply the procedure
described in \cite{GGG} to this new period matrix to find a
hyperelliptic curve $C$ over $K$ with
$(\Jac(C),\Theta)\stackrel{K}\simeq (A_f,\cL^{(\gamma^{-1})})$,
where $\Theta$ denotes the canonical principal polarization on $\Jac
(C)$. We obtain $Y^2=F(X)$,  where:
$$
F(X)=\frac{4\left(3-2\sqrt{-3}\right)}9X^6
+\frac{8(-1+\sqrt{-3})}3X^5 +\frac{4(3-7\sqrt{-3})}9X^4
+\frac{2(7+23\sqrt{-3})}{27}X^3$$
 $$
-\frac{(11+7\sqrt{-3})}{18}X^2 +\frac{(15+\sqrt{-3})}{108}X
+\frac{(-2+\sqrt{-3})}{324}.
$$
As we know that there is only one isomorphism class of principal
polarizations on $A_f$, by Torelli's theorem, the curves $C$ and
$\overline{C}$ must be isomorphic. Thus, their common Igusa
invariants must lie in $\Q$. Indeed, they are
$$
 \{i_1,i_2,i_3\}=\{\frac{2^{18}\cdot41^5}{3^3},  2^{12} \cdot3\cdot 41^{3},  2^9\cdot 7\cdot 41^2 \cdot47\}.
 $$

The curve $C$ has been obtained up to numerical approximations.
Next, we shall prove that $\Jac(C)$ is $K$-isogenous to the abelian
surface $A_f$. We point out  that the procedure which we will use to
determine the order $\End(\Jac(C))$ can be applied when $\Jac (C)$
has a unique class of principal polarizations and its endomorphism
algebra has a quaternion of reduced norm equal to $\pm 2$.

\begin{proposition} Let $C:Y^2=F(X)$ be the genus two curve given above.
Then, $\Jac(C)$ has multiplication by   a maximal order of the
quaternion algebra $\left ( \frac{6, -3}{\Q}\right)$ and is
$K$-isogenous to the abelian variety $A_f$.
\end{proposition}

\noindent {\bf Proof.} We split the proof in five steps.

\vskip 0.1 cm \noindent {\bf (i)} {\it  $C/K$ does  not admit a
descent to $\Q$}. It can be easily checked that $\Aut (C)=\langle
w\rangle$, where $w$ denotes the hyperelliptic involution. Let
$\overline C: \overline Y=\overline f (\overline X)$, where
$\overline f$ is the complex conjugate of $f$. Then, the isomorphism
$\varphi: C\longrightarrow \overline C$ given by
\begin{equation}\label{iso}  (\overline X, \overline Y)=\left(\frac{-1 + \sqrt{-3}
+ (3 - \sqrt{-3})X}{ 2\,\sqrt{-3} + 6(1 -
\sqrt{-3})X}\,,\,\frac{24\sqrt{-3}\,Y}{( 2\,\sqrt{-3} + 6(1 -
\sqrt{-3})X)^3}\right)
\end{equation}
satisfies that $\overline \varphi\circ \varphi=w$. By Weil's
criterion, this implies  that $C$ does
 not admits a descent over $\Q$.

\vskip 0.1 cm \noindent {\bf (ii)} {\it  $K$-linear relations
between the entries of a  period matrix  of $\Jac(C)$}.
 We will prove that a period matrix for the curve $C$ is
\begin{equation}\label{matriuperiodes}
  \Omega=\left (\begin{array}{cccc}
  \alpha& \displaystyle{-\frac{1+\sqrt{-3}}{2}\, \alpha - (3 - \sqrt{-3}) \beta}&\displaystyle{ \frac{-1+\sqrt{-3}}{2}\,\alpha +
  \frac{3-\sqrt{-3}}{2}\, \beta}& \displaystyle{\frac{3+\sqrt{-3}}{2}\, \beta}\\[5 pt]
  \beta &
 \displaystyle{-\frac{3+\sqrt{-3}}{6}\, \alpha - \frac{1 - \sqrt{-3}}{2} \beta}&
  \displaystyle{\frac{3+\sqrt{-3}}{12}\, \alpha - \frac{1 + \sqrt{-3}}{2} \beta}&
   \displaystyle{\frac{3 - \sqrt{-3}}{12}\,\alpha}
\end{array}\right)\,,
\end{equation}
where $\left(\alpha, \beta\right)=\int_{x_3}^{x_1}\left(\frac{dX}Y,
\frac{X\,dX}Y\right)$.

First of all, we remark that a symplectic basis for $H_1(C,\Z)$ is
given by the  pathes $\gamma_1$, $\gamma_2$, $\gamma_3$, $\gamma_4$,
enclosing respectively the segments $[x_1, x_3]$, $[x_4, x_5]$,
$[x_1, x_2]$, $[x_5, x_6]$, where
$$ \begin{array}{lll}
x_1=  0.085\cdots- 0.130 \cdots \,I\,, &x_2=0.098\cdots - 0.181
\cdots
\,I\,,&x_3=0.146\cdots - 0.232 \cdots \,I\,,\\
x_4=0.718\dots + 0.253 \dots \,I \,,& x_5=  0.751\dots - 0.187\dots
\,I \,,&x_6=
  0.770\dots - 0.018 \dots \,I\,,
\end{array}
$$
 are the roots of the polynomial $f(X)$ defining the curve $C$. It is well known how to build a period matrix for $C$. For instance, we can take:
 $$
\Omega=\left (\begin{array}{llll} {
 \int_{x_3}^{x_1}\frac{dX}{Y}}&{\int_{x_4}^{x_5}\frac{dX}{Y}}&
 {\int_{x_1}^{x_2}\frac{dX}{Y}}&{\int_{x_5}^{x_6}\frac{dX}{Y}}\\[6
 pt]
 {
 \int_{x_3}^{x_1}\frac{X\,dX}{Y}}&{\int_{x_4}^{x_5}\frac{X\,dX}{Y}}&
 {\int_{x_1}^{x_2}\frac{X\,dX}{Y}}&{\int_{x_5}^{x_6}\frac{X\,dX}{Y}}
\end{array}\right)\,,
$$
where we integrate along the segments joining the different roots of
$f(X)$, and the orientation of the pathes and the determination of $Y=\sqrt{f(X)}$ is taken  so that
 $$
 \Omega=
\left(
\begin{matrix}
 35.97...-7.80... I & -22.45...& -12.37...+21.43... I & 11.23... -7.80...I \\
 3.36... -7.14...I & - 15.73... & 2.25 + 3.90...I & 7.87... -
 7.14...I
\end{matrix}
\right) .
 $$

Three ingredients will be used to find relations between the entries
of this matrix:
\begin{enumerate}
\item[a)] The isomorphism $\varphi:
C\longrightarrow \overline C$ given in (\ref{iso}).

\item[b)] {\bf A Richelot isogeny}: If $\Jac (C)$ is simple over
$\overline \Q$, then for every subgroup $G$ of its $2$-torsion
isomorphic to $(\Z/2\Z)^2$, there is a genus two curve
$C':{Y'}^2=G(X')$ and an isogeny
$\phi:\Jac(C)\longrightarrow\Jac(C')$ whose kernel is $G$
and $\phi^*(X'^i\,dX'/Y')=X^i\,dX/Y$ for $0\le i\leq 1$. It  is called
the Richelot isogeny attached to $G$ (cf. Section 3.1 in
\cite{bostmestre}) or chapter 9 in \cite{CasselsFlynn}).

\item[c)] We know (cf. pag 90 in \cite{vig80}) that in all
maximal orders of  an indefinite quaternion algebra every integer
$d\in\Z$ is a reduced norm. In particular, if $\End (\Jac (C))$ is
such an order, by taking $d=2$ there is an endomorphism of $
\Jac(C)$ whose kernel $G$ is isomorphic to $(\Z/2\Z)^2$. Therefore,
if $\Jac(C)$ has a unique class of principal polarization then $C$
must be isomorphic to the curve $C'$ given by the Richelot isogeny
attached to $G$.
\end{enumerate}

In our case, to build a convenient Richelot isogeny, we take,
following Bost and Mestre  in \cite{bostmestre}:
$$
\begin{array}{l}
(P, Q, R):=\left(\frac{4\left(3-2\sqrt{-3}\right)}9(X-x_1)(X-x_3),
(X-x_2)(X-x_6),(X-x_4)(X-x_5)\right)
\end{array}
$$
and $\Delta$ equal to the determinant of the polynomials $P,Q,R$
with respect to the basis $1,X,X^2$, which turns to be $\Delta=
1/6$. We thus arrive at the genus two curve $C':Y'^2=1/\Delta U
(X')\,V(X')\, W(X')$ where $U=[Q,R]$, $V=[R,P]$ and $W=[P,Q]$, whose
Jacobian is isogenous to $\Jac C$. More precisely, $C'$ is given by
the equation $Y'^2=G(X')$, with:
$$
 G(X')= \frac{4(2 - \sqrt{-3})}{3} X'^6+ \frac{2(-9 + 7\,\sqrt{-3})}{3}X'^5+\frac{(2(5 - 9\,\sqrt{-3})}{3}X'^4
  +
$$
$$
 \frac{4(7 + 23\,\sqrt{-3})}{27}X'^3 + \frac{4(-9 -
5\,\sqrt{-3})}{27}X'^2+\frac{8}{27} X'+\frac{-3 +
2\,\sqrt{-3}}{243}\,.
$$
Let
$$\begin{array}{ll}
x'_1=0.040\cdots - 0.214 \cdots\,I\,,&x'_2= 0.112\cdots - 0.184 \cdots\,I\,,\\
x'_3=0.149\cdots - 0.154 \cdots\,I\,, &x'_4= 0.606\cdots - 0.024 \cdots\,I\,,\\
x'_5=0.808\cdots + 0.019\cdots\,I\,, &x'_6=1.069\cdots - 0.059 \cdots\,I\,,
\end{array}
$$
be the roots of the polynomial $G(X')$ defining the curve $C'$. The
description of the  isogeny $\Jac C\rightarrow \Jac C'$ in terms of
the periods of $C$ and $C'$ is given by the equalities ($i=0,1$):
\begin{equation}
\label{igualtats-isogenia}
\begin{array}{lllllll}
\ds\int_{x_1}^{x_3} \frac{X^i\,dX}{Y}&= &\ds2 \int_{x'_2}^{x'_3}
\frac{X'^i\,dX'}{Y'}&\phantom{ccccc}&\ds\int_{x_1}^{x_2}
\frac{X^i\,dX}{Y}&= &
\ds\int_{x'_3}^{x'_1}    \frac{X'^i\,dX'}{Y'}\\[3pt]
\ds\int_{x_2}^{x_6} \frac{X^i\,dX}{Y}&= &\ds2 \int_{x'_3}^{x'_4}
\frac{X'^i\,dX'}{Y'}& &\ds\int_{x_3}^{x_5} \frac{X^i\,dX}{Y}&= &
\ds\int_{x'_2}^{x'_5}
\frac{X'^i\,dX'}{Y'}\\[3pt]
\ds\int_{x_4}^{x_5} \frac{X^i\,dX}{Y}&= &\ds2 \int_{x'_5}^{x'_4}
\frac{X'^i\,dX'}{Y'}& &\ds\int_{x_5}^{x_6} \frac{X^i\,dX}{Y}&= &
\ds\int_{x'_4}^{x'_6} \frac{X'^i\,dX'}{Y'}
\end{array}
\end{equation}
The trick to prove these equalities is the following: by
construction, six of the integrals $\int_{x_j}^{x_k} \frac{X^idX}{Y}$ on the left must
be equal to either an integral $\int_{x'_{j'}}^{x'_{k'}}
\frac{X'^idX'}{Y'}$ or  $2\,\int_{x'_{j'}}^{x'_{k'}}
\frac{X'^idX'}{Y'}$ for some $j', k'$; since these are all of them different, only a match is possible between the integrals on $C$ and the integrals on $C'$: evaluating  all the integrals to enough accuracy, we discard all but one possible equality, which must be true.

We can now relate the integrals  $\int_{x'_{j'}}^{x'_{k'}}
\frac{X'^idX'}{Y'}$ with the periods of $C$  by using the
isomorphism $\phi :C \longrightarrow C'$ given by
$$
(X,Y)\mapsto (X',Y')=\left(\frac{1 - \sqrt{-3}}{12\,X}\,, \quad
\frac{\sqrt{-3} \,Y}{18\, X^3}\right).
$$
We have
$$
\phi(x_1)=x'_6,\quad \phi(x_2)=x'_5, \quad \phi(x_3)=x'_4, \quad
\phi(x_4)=x'_1,\quad \phi(x_5)=x'_3, \quad\phi(x_6)=x'_2,
$$
and
$$
\phi^\ast\left(\frac{dX'}{Y'}\right)=\frac{3+\sqrt{-3}}{2}\frac{X\,dX}{Y},\quad
\phi^\ast\left(\frac{X'dX'}{Y'}\right)=\frac{3-\sqrt{-3}}{12}\frac{dX}{Y}.
$$
From these relations we obtain, for instance:
$$
\int_{x'_3}^{x'_1}\frac{dX'}{Y'}=\frac{3+\sqrt{-3}}{2}\int_{\phi^{-1}[x'_1,x'_3]}
\frac{X\,dX}{Y}.
$$
The transformed path $\phi^{-1}[x'_3,x'_1]$ is drawn on the
following graphic:
\begin{center}
\includegraphics[height=4cm]{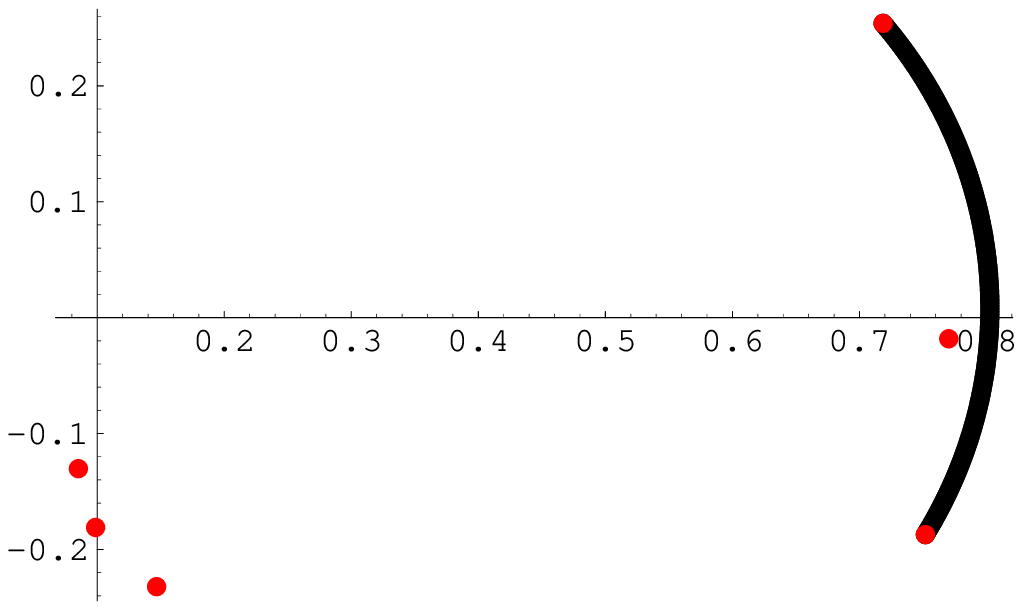}
\end{center}
Hence,  it is homologically equivalent to the path $[x_4, x_6]+[x_6,
x_5]$, so that
$$
\int_{x'_1}^{x'_3}\frac{dX'}{Y'}=\frac{3+\sqrt{-3}}{2}\left(\int_{x_4}^{x_6}+\int_{x_6}^{x_5}
\right)\frac{X\,dX}{Y}.
$$
Taking into account relations (\ref{igualtats-isogenia}), we obtain:
$$
\int_{x_1}^{x_2}\frac{dX}{Y}=\frac{3+\sqrt{-3}}{2}\left(\int_{x_4}^{x_6}+\int_{x_6}^{x_5}
\right)\frac{X\,dX}{Y}.
$$
In a similar way one can prove the following equalities:
$$
\begin{array}{ll}
\ds\int_{x_1}^{x_2}\frac{dX}{Y}=\frac{3+\sqrt{-3}}{2}\left(\int_{x_4}^{x_6}+\int_{x_6}^{x_5}
\right)\frac{X\,dX}{Y},
&\ds\int_{x_1}^{x_2}\frac{X\,dX}{Y}=\frac{3-\sqrt{-3}}{12}\left(\int_{x_4}^{x_6}+\int_{x_6}^{x_5} \right)\frac{dX}{Y},\\\\
\ds\int_{x_4}^{x_5}\frac{dX}{Y}=-\frac{3+\sqrt{-3}}{2}\int_{x_2}^{x_3}\frac{X\,dX}{Y},
&\ds\int_{x_4}^{x_5}\frac{X\,dX}{Y}=-\frac{3+\sqrt{-3}}{12}\int_{x_2}^{x_3}\frac{dX}{Y},\\\\
\ds\int_{x_5}^{x_6}\frac{dX}{Y}=-\frac{3+\sqrt{-3}}{2}\int_{x_1}^{x_3}\frac{X\,dX}{Y},
&\ds\int_{x_5}^{x_6}\frac{X\,dX}{Y}=-\frac{3+\sqrt{-3}}{12}\int_{x_1}^{x_3}\frac{dX}{Y}.
\end{array}
$$
The last line gives immediately the expression of the last column of
the period matrix $\Omega$ in terms of the first column. To relate
the third column with the first and second  columns of $\Omega$, we
substitute the equalities
$\int_{x_2}^{x_3}=\int_{x_2}^{x_1}-\int_{x_1}^{x_3}$ and
$\int_{x_4}^{x_6}=\int_{x_4}^{x_5}+\int_{x_5}^{x_6}$ in the first
two equalities above, and look at the resulting equalities as a
linear system of equations with respect to the integrals
$\int_{x_1}^{x_2}$. Solving it, we find:
$$
\begin{array}{l}
\ds\int_{x_1}^{x_2}\frac{dX}{Y}=\frac{3+\sqrt{-3}}2\int_{x_4}^{x_5}\frac{X\,dX}{Y}+\int_{x_1}^{x_3}\frac{X\,dX}{Y},\\\\
\ds\int_{x_1}^{x_2}\frac{X\,dX}{Y}=\frac{-3+\sqrt{-3}}{12}\int_{x_4}^{x_5}\frac{dX}{Y}+\int_{x_1}^{x_3}\frac{dX}{Y}.\\\\
\end{array}
$$
Hence, we have proved that
\begin{equation}\label{eqs1}
  \Omega=\left (\begin{array}{cccc}
  \alpha& \delta&\displaystyle{ -\beta+\frac{3+\sqrt{-3}}{2}\,\gamma
  }& \displaystyle{\frac{3+\sqrt{-3}}{2}\, \beta}\\[5 pt]
  \beta &
 \gamma&
  \displaystyle{-\alpha+\frac{-3+\sqrt{-3}}{12}\, \delta}&
   \displaystyle{\frac{3 - \sqrt{-3}}{12}\,\alpha}
\end{array}\right)\,,
\end{equation}
 for some $\alpha,\beta,\delta,\gamma\in\C$.  Unfortunately, we do not obtain more information when we make a similar treatment from the intermediate equalities.
At this point, plays its role the first ingredient of our proof, the
isomorphism $\varphi: C\rightarrow \overline{C}$. Note that
 \begin{equation}\label{integral-conjugada}
 \overline{\int_{x_4}^{x_5}\frac{X^i dX}Y}=\int_{\overline{x_4}}^{\overline{x_5}}\frac{\overline{X} ^i d\overline{X} }{\overline{Y}},
 \end{equation}
 and we can pull back the last integrals to periods of $C$ by means of
 $\varphi$. The path $\varphi^{-1}[\overline{x_4},\overline{x_5}]$ is
homologous to $-[x_1,x_2]+[x_2,x_3]=-2\,[x_1,x_2]+[x_1,x_3]$, as
seen in the following graphic:
\begin{center}
\includegraphics[height=4cm]{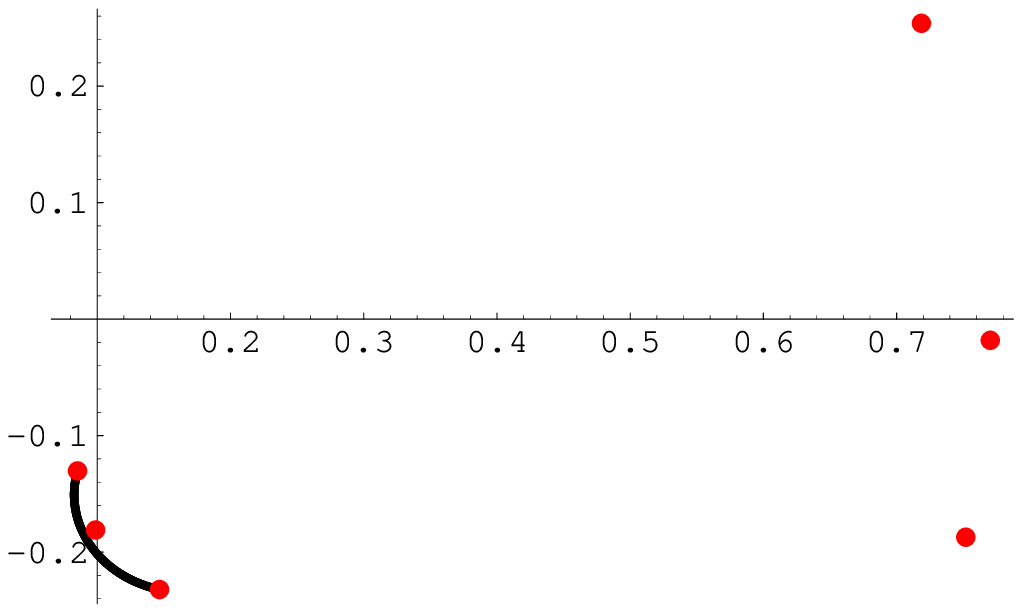}
\end{center} We obtain:
$$
\int_{\overline{x_4}}^{\overline{x_5}}\frac{\overline{X} ^i
d\overline{X} }{\overline{Y}}
=\int_{\varphi^{-1}[\overline{x_4},\overline{x_5}]}\varphi^\ast
\frac{\overline{X} ^i d\overline{X} }{\overline{Y}}=
\left(-2\int_{x_1}^{x_2}+\int_{x_1}^{x_3}\right)\varphi^{\ast}\frac{\overline{X}
^i d\overline{X} }{\overline{Y}},
$$
where
$$
\varphi^\ast\frac{ d\overline{X}
}{\overline{Y}}=-\frac{1+\sqrt{-3}}2\frac{dX}Y+2\sqrt{-3}\frac{X\,dX}Y,\quad\
\varphi^\ast\frac{\overline{X} d\overline{X}
}{\overline{Y}}=\frac{-\sqrt{-3}}3\frac{dX}Y-\frac{1-\sqrt{-3}}2\frac{X\,dX}Y.
$$
Substituting these equalities in (\ref{integral-conjugada}) and
equating real and imaginary parts, we obtain four independent linear
equations in the variables $\re(\alpha)$, $\im (\alpha)$, $\re(\beta)$,
$\im (\beta)$, $\re(\delta)$, $\im (\delta)$, $\re(\gamma)$ and $ \im
(\gamma)$. The solution of  this system depending on the variables
$\re(\alpha)$,
 $\im (\alpha)$, $\re(\beta)$ and $ \im (\beta)$ yields the equality
(\ref{matriuperiodes}).

\vskip 0.1 cm \noindent {\bf (iii)} {\it The ring $\End_K(\Jac(C))$
is a maximal order of $\left ( \frac{6, -3}{\Q}\right)$}.   It can be checked now that the matrices:
$$
T_6= \left(\begin{array}{cc} 0 & 1\\ 6& 0\end{array}\right)\,,\quad
T_{-3}= \left(\begin{array}{cc} \sqrt{-3} & 0\\0&
-\sqrt{-3}\end{array}\right)
$$
satisfy $T_6^2=6$, $T_{-3}^2=-3$, $T_6\cdot T_{-3}=-T_{-3}\cdot T_6$
and, , by
using (\ref{matriuperiodes}), we have that:
$$
\Omega \cdot
\left(\begin{array}{rrrr}
-2& -2& 1& 0\\ -2& 2& -1& -1\\
 -2 & 0& 0 & -2\\-2 & 2 & -4 &0
\end{array}\right)={}^tT_6 \cdot \Omega\,,
\quad\quad \Omega\cdot
\left(\begin{array}{rrrr}
3 & -2 & -2& 2\\ 2& -3& 0& 2\\
 4 & 4 & -1&2\\0 & -4 & 2& 1
 \end{array}\right)={}^tT_{-3}\cdot \Omega\,,
$$
which implies that there are $i, j\in \End_K(\Jac (C))$ such that
the matrices of their actions on $H^0(C,\Omega^1)$ with respect the
basis $\{dX/Y, X\,dX/Y\}$ are $T_6$ and $T_{-3}$ respectively (see
Lemma \ref{criterion}). Moreover, it is easy to check that  the
maximal order $\langle 1,\frac{1}{2}i+\frac{1}{6}
k,\frac{1}{2}+\frac{1}{2} j, \frac{1}{3} k\rangle $ of $\Q(i,j)$ is
contained in $\End_K(\Jac (C))$. To assert that this maximal order
is the full ring endomorphisms of $\Jac (C)$, we only need to prove
that $C$ is not isogenous to a square of an elliptic curve with complex multiplication .
But, it follows from the fact that the reduction of $C$ at a prime
of $K$ dividing to $7$ (resp. $13$) is isogenous to the square of an
ordinary elliptic curve with multiplication by $\Q(\sqrt{-6})$
(resp. $\Q(\sqrt{-51})$).

\vskip 0.1 cm \noindent {\bf (iv)} {\it  $\Jac(C)$  is
$K$-isomorphic to an abelian surface $A/\Q$ with $\End_\Q(A)\simeq
\Z[\sqrt 6 ]$.} Let $u=1/2- j/2+k/3\in\Aut (\Jac(C))$. Let $\iota$
be the isomorphism obtained as the composition of the following
isomorphisms
$$\Jac(C)\stackrel{u}\longrightarrow \Jac(C)\stackrel{\varphi_*}\longrightarrow\Jac(\overline C),$$
where  $\varphi$ is as in (i). The matrix of the
morphism
$$\iota^*:H^0( \Jac(\overline C),\Omega^1_{\Jac(\overline
C)/K})\longrightarrow H^0( \Jac(C),\Omega^1_{\Jac(C)/K})$$ with
respect the bases $\{d \overline X/\overline Y,\overline X\,d
\overline X/\overline Y\}$ and $\{dX/Y,X\,dX/Y\}$ is the identity
matrix.
Therefore, $\overline{\iota}\circ\iota=1$ and, then, by the Weil
criterion there is an abelian surface $A/\Q$ and an isomorphism
$\mu:\Jac (C)\longrightarrow A$  defined over $K$ such that
$\overline \mu=\mu\circ\iota^{-1}$. Since $\nu=\mu\circ
i\circ\mu^{-1}$ is invariant under the complex conjugation, it must
be  $\Z [\nu]\subseteq \End_\Q(A)$. Due to the fact that $\nu^2=6$
and $\Z[\sqrt{6}]$ is a maximal order, it follows that $\Z [\nu]=
\End_\Q(A)$.

\vskip 0.1 cm \noindent {\bf (v)} {\it $A$ is $\Q$-isogenous to
$A_f$.} The conductor of $\Jac (C)$ is the integer ideal of $K$
generated by $3^{8}$. By the recent progress with respect to the
generalized Shimura-Taniyama-Weil Conjecture (cf. \cite{KW06}),  $A$ is
modular and, thus, there exists a newform $g$ of level $M$ such that
$A$ is $\Q$-isogenous to $A_g$. Since $g$ must have an inner-twist
by the Dirichlet character $\chi$ attached to $K$,i.e. $g\otimes
\chi={}^\sigma g$, we have:
$$
\Res_{K/\Q}(\Jac (C))=\Res_{K/\Q}(A)\stackrel{\Q}\sim A_g\times
A_{g\otimes \chi}=A_g^2\,.
$$
 By applying Milne's formula for the conductor of
the Weil restriction (cf. Proposition 1 in \cite{milne}) and the
fact that $\operatorname{cond\,}(A_g)=M^2$ (cf. \cite{carayol}), we
obtain that $\operatorname{cond\,}(A_g)=3^{10} $. The assertion is
obtained by checking that in the level $243=3^5$ the unique
normalized newforms $g$ such that $A_g$ is a quaternionic  surface
are $f$ and its Galois conjugate. \hfill $\Box$

\begin{rem}
Since $2$ divides the discriminant of the quaternion algebra $\End^0
(\Jac(C))$, there is a unique ideal of norm $2$ in any of its
maximal orders. So, all quaternions of reduced norm $2$, viewed as
endomorphisms of $\Jac(C)$, have the same kernel and, thus, only one
of the $15$ curves obtained through a Richelot-isogeny is isomorphic
to $C$.
\end{rem}


\section{Example: $A_f$ with no principal polarizations}

We can apply the same process to the normalized newform in $S_2(\Gamma_0(972))$:
$$
f=
q + 3\sqrt{2} q^5 + 2q^7 + 3\sqrt{2} q^{11} - q^{13} - 3\sqrt{2} q^{17} + 5q^{19}
 + \cdots
$$
The  canonical polarization  on $J_0(N)$ induces a polarization
$\cL$ on $A_f$ of type $[1,6]$. The endomorphism algebra of $A_f$ is
$\End_K(A_f)=\langle  1,\frac{3}{2}i+\frac{1}{2}
k,\frac{1}{2}+\frac{1}{2} j,  k\rangle$, with $i^2=2, j^2=-3$. It
has index 3 with respect to  a maximal order in $\End^0_K(A_f)$.
Since there are no symmetric totally positive elements of reduced
norm 6, we conclude that $A_f$ admits no principal polarizations.

\section{Example:  $A_f$ with two principal polarizations}

Let us finally consider the abelian surface $A_f$ attached to the normalized newform in $S_2(\Gamma_0(1568))$
$$
f=q + \sqrt{7}q^3 - 3q^5 + 4q^9 - \sqrt{7}q^{11} - 4q^{13} - 3\sqrt{7}q^{15} + \cdots
$$
which corresponds to the entry labelled $S_{1568L}$ in  Table II.
Now $i^2=7, j^2=-1$, and we find that
 $\End_K(A_f)=\Z\langle1,
\frac12+\frac12i+\frac12j+\frac12k, j,k\rangle$ is a maximal order
in $\Q(i,j)$. The polarization $\cL$  on $A_f$ induced by the
canonical polarization on $J_0(1568)$ is of type $[1,14]$.  By
theorem \ref{rotger}, there are two isomorphism classes of principal
polarizations on $A_f$ over $K=\Q(\sqrt{-1})$, and by lemma
\ref{pp}, they are not defined over $\Q$. Thus, there exist two
$C_1, C_2$ over $K$ with $\Jac(C_i)\simeq A_f$. They correspond to
the principal polarizations $\cL^{(\gamma_1^{-1})}$ and
$\cL^{(\gamma_2^{-1})}$, where
$$\gamma_1=7+2\,i+k \quad \text{and}\quad  \gamma_2=7-2\,i+k\,.$$
A priori, both $C_1$ and
$C_2$ could be $K$-isomorphic to their respective complex conjugates
as above; if this is not the case, then they must be complex
conjugates. Indeed, computations reveal that we are in this second
case: the Igusa invariants of $C_1$ and $C_2$ are
$$
\left(
\begin{array}{c}
i_1(C_1)\\i_2(C_2)\\i_3(C_3)
\end{array}\right)=
\left(\begin{array}{c} \frac{  ( 1 + I)^{14}( -7 + 8I)^5 ( 28+
5I)^5 }{( 2 + I)^{12}} \\
 \frac{  ( 1 + I)^{10} ( 3 +10I)^2 (7-8 I)^3(28+5I)^3} {( 2 +
 I)^8}\\
 \frac {( 1 + I)^{12}(-2+3I)(8+7I)^2(28+5I)^2(320+1383I)} {( 2+ I)^8}
\end{array}\right)=
\left(
\begin{array}{c}
\overline{i_1(\overline{C_1})}\\\overline{i_2(\overline{C_2})}\\\overline{i_3(\overline{C_3})}
\end{array}\right).
$$

A rational equation for $C_1$ is:
$$
\begin{array}{rl}
Y^2=\ds\frac{1372-539 I}{5}& \ds \left(X^6+\frac{332+208 I}{181}
   X^5+\frac{1173+2148 I}{1267} X^4+\frac{376+8060 I}{8869}
   X^3\right.\\\\
&\ds \left.  -\frac{705-1992 I}{8869} X^2-\frac{1228-1612 I}{62083}
   X-\frac{607-492 I}{434581}\right)
\end{array}
$$

It is worth mentioning that up to level $N\le 2500$,  there are only
two modular abelian surfaces  $A_{f_1}$ and $A_{f_2}$
 with quaternionic multiplication
and more than one isomorphism class of principal polarizations. Both
newforms have level $N=1568$ and, moreover, each of them is the
twist of the other  by the quadratic Dirichlet character of
conductor $28$.

\section{Table of quaternionic abelian surfaces}

Using a program made by J. Quer for {\em Magma} we have found  all
quaternionic abelian surfaces $A_f$ which appear as subvariety of
some $J_0(N)$ with $N\leq 7000 $. In the first table below, we list
for every level $N$ the number of factors $A_f$ with attached number
fields $E, K$. The reduced discriminant of $\End^0(A_f)$ is denoted
by~$D$.

$$
\begin{array}{c}
\begin{array}{|ccccc|}
N & E & K & D& \# A_f\\
\hline \hline
243 & \Q( \sqrt {6}) & \Q( \sqrt{-3})&6& 1\\[3 pt]
675 & \Q( \sqrt {2}) & \Q( \sqrt{-3})& 6&2\\[3 pt]
972 & \Q( \sqrt {2}) & \Q( \sqrt{-3})&6& 1\\[3 pt]
1323 & \Q( \sqrt {6}) & \Q( \sqrt{-3})& 6&2\\[3 pt]
1568 & \Q( \sqrt {7}) & \Q( \sqrt{-1})&14& 2\\[3 pt]
1568 & \Q( \sqrt {3}) & \Q( \sqrt{-1})& 6&2\\[3 pt]
1849 & \Q( \sqrt {6}) & \Q( \sqrt{-43})&6& 1\\[3 pt]
2592 & \Q( \sqrt {6}) & \Q( \sqrt{-1})& 6&2\\[3 pt]
2592 & \Q( \sqrt {3}) & \Q( \sqrt{-1})&6& 2\\[3 pt]
2601 & \Q( \sqrt {2}) & \Q( \sqrt{-51})& 6&1\\[3 pt]
2700 & \Q( \sqrt {10}) & \Q( \sqrt{-3})& 10&2\\[3 pt]
3136 & \Q( \sqrt {3}) & \Q( \sqrt{-1})&6& 2\\[3 pt]
 \hline
\end{array}
\quad
\begin{array}{|ccccc|}
N & E & K &  D&\# A_f\\
\hline \hline
3136 & \Q( \sqrt {7}) & \Q( \sqrt{-1})& 14&2\\[3 pt]
3886 & \Q( \sqrt {6}) & \Q( \sqrt{-3})& 6&1\\[3 pt]
3886 & \Q( \sqrt {2}) & \Q( \sqrt{-3})& 6&1\\[3 pt]
3969 & \Q( \sqrt {15}) & \Q( \sqrt{-7})& 15&1\\[3 pt]
5184 &  \Q(\sqrt 3) & \Q(\sqrt{-1})& 6& 2\\[3pt]
5184 &  \Q(\sqrt 6) & \Q(\sqrt{-3})& 6&2\\[3pt]
5292 & \Q(\sqrt{10})&\Q(\sqrt{-3}) & 10&2\\[3pt]
5408 & \Q(\sqrt{11})&\Q(\sqrt{-1}) & 22&2\\[3pt]
5408 & \Q(\sqrt{3})&\Q(\sqrt{-13}) & 6&2\\[3pt]
6075 & \Q( \sqrt {6}) & \Q( \sqrt{-3})&6& 2\\[3 pt]
6400 & \Q( \sqrt {6}) & \Q( \sqrt{-1})& 6&4\\[3 pt]
   &       &     & &\\[3 pt]
\hline
\end{array}\\\\
\mbox{ Table I}
\end{array}
$$

For the $A_f$  in the table above with level $N\le2500$, we describe its endomorphism algebra and its principal polarizations in the following table.
We follow the
labelling of {\em Magma} to denote the $A_f$; for instance
$S_{243,D}$ denotes the fourth $\Q$-irreducible factor of $J_0(243)$
For every $A_f$ we describe the quadratic fields $E=\Q(\sqrt{m})$,
$K=\Q(\sqrt{\delta})$, the endomorphism algebra $\End_K(A_f)$, its
index $n$ with respect to a maximal order in $\End_K^0(A_f)$ and the
numbers $\pi_\Q$,  $\pi_K$ of polarizations over $\Q$ and $K$
respectively. For those $A_f$ admitting principal polarizations, we
provide elements $\gamma$ originating them.

$$\begin{array}{cccccccc}
A_f&(m,\delta) & d & \End_K(A_f) & n &\pi_\Q&\pi_K&\gamma\\[4,pt]
\hline \hline S_{243D} & (6,-3) & 2 &\langle 1,
\frac{1}{2}i+\frac{1}{6} k,\frac{1}{2}+\frac{1}{2} j, \frac{1}{3} k
\rangle      & 1& 0&
1&2+\frac12 i-\frac16k\\[4 pt]
\hline S_{675L} & (2,-3) & 2 & \langle 1, \frac{1}{2}i+\frac{1}{2}
k,\frac{1}{2}+\frac{1}{2} j,  k \rangle      & 1& 1&
1 &2+i\\[4 pt]
\hline S_{675P} & (2,-3) & 2 & \langle  1, \frac{1}{2}i+\frac{1}{2}
k,\frac{1}{2}+\frac{1}{2} j,  k\rangle      & 1& 1&
1 &2+i\\[4 pt]
\hline S_{972E} & (2,-3) & 6 & \langle  1,\frac{3}{2}i+\frac{1}{2}
k,\frac{1}{2}+\frac{1}{2} j,  k\rangle      & 3& 0& 0 &\\[4 pt]
\hline S_{1323U} & (6,-3) & 6 & \langle 1,  \frac{1}{2}i+\frac{1}{2}
k,\frac{1}{2}+\frac{1}{2} j,  k\rangle      & 3& 0&
0 &\\[4 pt]
\hline S_{1323V} & (6,-3) & 6 & \langle  1, \frac{1}{2}i+\frac{1}{2}
k,\frac{1}{2}+\frac{1}{2} j,  k\rangle
& 3& 0&0 &\\[4 pt]
\hline S_{1568L} & (7,-1) & 14 & \langle  1,
\frac12+\frac12i+\frac12j+\frac12k, j,k\rangle & 1& 0& 2
&\begin{array}{c} 7+2i+k\\7-2i+k
\end{array}\\[4 pt]
\hline S_{1568N} & (7,-1) & 14& \langle 1,
\frac12+\frac12i+\frac12j+\frac12k, j,k \rangle      & 1& 0& 2
&\begin{array}{c} 7+2i+k\\7-2i+k
\end{array}\\[4 pt]
\hline S_{1568S} & (3,-1) & 6 & \langle 1,
\frac12+\frac12i+\frac12j+\frac12k, j,k  \rangle      & 1& 0&
1 &3+k\\[4 pt]
\hline S_{1568U} & (3,-1) & 6 & \langle 1,
\frac12+\frac12i+\frac12j+\frac12k, j,k  \rangle      & 1& 0&
1 &3+k\\[4 pt]
\hline S_{1849E} & (6,-43) & 258 & \langle
1,\frac12i+\frac12k,\frac12+\frac12j,k \rangle & 43& 0&
0 &\\[4 pt]
\hline\\[1pt]
&&&\mbox{ Table II}
\end{array}
$$

\providecommand{\bysame}{\leavevmode\hbox
to3em{\hrulefill}\thinspace}
\providecommand{\MR}{\relax\ifhmode\unskip\space\fi MR }
\providecommand{\MRhref}[2]{%
  \href{http://www.ams.org/mathscinet-getitem?mr=#1}{#2}
}

\end{document}